\newcommand{\bull}{\vrule height .9ex width .8ex depth -.1ex}
\newcommand{\ppp}{\hfill $\bull$ }
 \documentclass[11pt]{article}
 \author{ Mohammed Larbi Labbi}
   \title{On Generalized  Einstein Metrics}
   \date{}
   \usepackage{amsmath}
   \usepackage{amscd}
\newtheorem{theorem}{Theorem}[section]
\newtheorem{corollary}[theorem]{Corollary}

\newtheorem{proposition}[theorem]{Proposition}
\newtheorem{definition}{Definition}[section]

\begin{document}
   \maketitle
   \begin{abstract} Recall that the usual Einstein metrics are those for which the first Ricci
 contraction of the covariant Riemann curvature tensor is proportional to the metric. 
Assuming the same type of restrictions 
but instead on the different contractions of
Thorpe tensors, one gets several natural generalizations 
of Einstein's condition. In this paper, we study some properties of these classes of metrics. 

\end{abstract}
\section{Introduction}
Let $(M,g)$ be a Riemannian manifold of dimension $n$ and $R$ be its Riemann covariant curvature tensor.
For each even integer $2p$,  $2\leq 2p\leq n$, J. A. Thorpe \cite{Thorpe} introduced a generalization of the 
Gauss-Kronecker
curvature tensor as follows: For each $m\in M$ and for $u_i,v_j$ tangent vectors at $m$ we set\\
\begin{equation*}
\begin{split}
R_{2p}(u_1,...,u_{2p},v_1,...,v_{2p})=\frac{2^{-p}}{(2p)!}
\sum_{\scriptstyle \alpha, \beta \in S_{2p}}
&\epsilon(\alpha)\epsilon(\beta)
 R(u_{\scriptstyle \alpha(1)},u_{\scriptstyle \alpha(2)},
v_{\scriptstyle \beta(1)},v_{\scriptstyle \beta(2)}) ...\\
&... R(u_{\scriptstyle \alpha(2p-1)},u_{\scriptstyle \alpha(2p)},
v_{\scriptstyle \beta(2p-1)},v_{\scriptstyle \beta(2p)}).
\end{split}
\end{equation*}
Where $S_{2p}$ denotes the group of permutations of $\{1,...,2p\}$ and, for $\alpha\in S_{2p}$,
$\epsilon(\alpha)$ is the signature of $\alpha$. In particular, $R_2$ is the Riemann covariant curvature tensor
$R$. If the dimension $n$ is even, $R_n(e_1,...,e_n,e_1,...,e_n)$ is the Gauss-Bonnet
integrand of $(M,g)$, $\{e_1,...,e_n\}$ being an orthonormal basis of $T_mM$.\\
It is evident that the tensor $R_{2p}$ is alternating in the first $(2p)$-variables, alternating
in the last $(2p)$ variables, and is invariant under the operation of interchanging the 
first $2p$
variables with the last $2p$. Hence it is a symmetric $(2p,2p)$ double form. Using the exterior
product of double forms, the tensors $R_{2p}$ can be written in the following compact 
form \cite{Kulk, Nasu, Labbidoubleforms},
\begin{equation*}
R_{2p}=\frac{2^{p}}{(2p)!}R^{p}.
\end{equation*}
Where $R^{p}$ denotes the exterior product of $R$ with itself $p$ times.\\
The full Ricci contraction of $R_{2p}$ is, up to a constant, the $(2p)$-th Gauss-Bonnet curvature of 
$(M,g)$ \cite{Labbivariation, Labbidoubleforms}.\\
For $0< k< 2p$, we shall say that a metric is   $(k,p)$-Einstein if the $k^{th}$ Ricci contraction 
of its Thorpe curvature tensor $R_{2p}$  is proportional to $g^{2p-k}$. This condition
generalizes the usual Einstein condition obtained for $p=k=1$, and it
implies that the metric is $2p$-Einstein in the sens of \cite{Labbivariation}.
 In particular, in the compact case the $(k,p)$-Einstein metrics are for any $k$   
critical metrics of  the total $2p$-th Gauss-Bonnet
curvature functional once restricted to metrics of unit volume \cite{Labbivariation}.\\
The paper is divided into three parts. In the first one, we recall some useful facts 
about double forms and by the way we  complete and improve two results of 
\cite{Labbidoubleforms}. These results are then used in the second part where we study some geometric properties of
$(k,p)$-Einstein metrics. We emphasize in the second part the special class of $(1,p)$-Einstein manifolds, these 
manifolds behave in many directions like the usual $4$-dimensional Einstein manifolds. In the last part,
we discuss and generalize a related but more subtle  generalization of Einstein metrics suggested by a work of
 Thorpe.
\section{Preliminaries} 
Let $(M,g)$ be a smooth Riemannian manifold of dimension n. A double form of degree $p$ 
is a section of the bundle ${\cal D}^{p,p}= \Lambda^{*p}M \otimes
   \Lambda^{*p}M$ where  
 $\Lambda^{*p}M$ denotes the bundle of differential
  forms of degree $p$ on $M$. Double forms are abundant in Riemannian geometry:
 the metric, Ricci
and Einstein tensors are symmetric double forms of degree 1, The Riemann curvature tensor, 
Weyl curvature tensor are symmetric double forms of order two, Gauss-Kronecker and 
Weitzenb\"ock curvature tensors 
are examples of symmetric double forms of higher order...\\
 The usual exterior product of differential 
forms extends in a natural way
way to double forms. The resulting product is some times called 
 Kulkarni-Nomizu product of double forms.
In particular, $\frac{k}{2}$ times the exterior product of the metric with itself  
$\frac{k}{2}g^2$ is the curvature tensor of a Riemannian 
manifold of constant curvature $k$, the repeated products of the Riemann curvature tensor
 $RR...R=R^q$ determine the generalized Gauss-Kronecker tensors defined above.  \\
If we denote by $E^{(r,r)}$ the trace free double forms of order $r$, then we have the 
following orthogonal decomposition of the bundle
$D^{(p,p)}$ for $1\leq p\leq n$:
\begin{equation}\label{ort:decom}
D^{p,p}=E^{p,p}\oplus gE^{p-1,p-1}\oplus g^2E^{p-2,p-2}\oplus ...
\oplus g^pE^{0,0}.
\end{equation}
The following proposition is a direct consequence of corollary 2.4 in \cite{Labbidoubleforms}.
\begin{proposition}\label{firstproposition} If $n<2p$, then 
\begin{equation}\label{firstequation}
D^{p,p}=g^{2p-n}D^{n-p,n-p}
\end{equation}
In particular, if $\omega=\sum_{i=0}^{p}g^i\omega_{p-i}\in D^{p,p}$ is the decomposition 
of a double form $\omega$ following the orthogonal splitting (\ref{ort:decom}), then
\begin{equation}\label{secondequation}
\omega_k=0\, \, \, {\rm for}\, \, \, n-p<k\leq p.
\end{equation}
\end{proposition}
{\sl Proof.} Straightforward. Just set $j=2p-n$ and discuss following the cases where $j$ is 
 even or  odd 
using corollary 2.4 of 
\cite{Labbidoubleforms}.\ppp
\par\medskip\noindent
{\sc Remark.} For a given double form $\omega$ of order $p$, we evaluated 
in theorem 3.7 of \cite{Labbidoubleforms}, the components
 $\omega_k$ with respect to the orthogonal splitting (\ref{ort:decom}). The author forgot
to include in the assumptions  of that  formula the condition $n\geq 2p$. This condition is
 necessary as if not 
some coefficients will vanish in the  proof. Also, it is clear that the coefficients
$(n-2k)!$ and $(n-p-k)!$ in the final formula are undefined if we do not assume $n\geq 2p$.\\
However, if $n<2p$, the previous proposition asserts that $\omega$ is divisible by
some power of the metric $g$ and hence some components vanish as above, the other components
can be determined by applying  theorem 3.7  to the quotient double form.
\par\medskip\noindent
The following proposition improves lemma 5.7 of \cite{Labbidoubleforms} and its proof below
completes a missing part in the proof of the same lemma.
\begin{proposition}\label{secondproposition} Let $\omega=\sum_{i=0}^{p}g^i\omega_{p-i}$ be a double form of order 
$p$ and $1\leq k<p$. Then $c^k\omega$
is proportional to the metric $g^{p-k}$ if and only if
\begin{equation*}
\omega_r=0\, \, \, \, {\rm for}\, \, \, \, 1\leq r\leq p-k.\end{equation*}
\end{proposition}
{\sl Proof.}   Formula (12) in \cite{Labbidoubleforms} shows that
$$c^k\omega=\sum_{i=k}^p\alpha_{ik}g^{i-k}\omega_{p-i},$$
where $\alpha_{i0}=1$ for $0\leq i\leq p$. For $1\leq k\leq i\leq p$ we have
$$\alpha_{ik}=\frac{i!}{(i-k)!}\prod_{j=1}^{j=k}(n-2p+i+j).$$
Remark that all the coefficients $\alpha_{ik}$ are nonnegative for $0\leq k\leq i\leq p$, furthermore,
$\alpha_{ik}=0$ if and only if $2p>n$ and $i<2p-n$. In such a case we do have $\omega_{p-i}=0$
by proposition \ref{firstproposition}.\\
Now, suppose $\omega_r=0$  for $1\leq r\leq p-k$, then $\omega_{p-i}=0$ for $k\leq i\leq p-1$ and
therefore $c^k\omega=\alpha_{ip}g^{p-k}\omega_0$ is proportional to the metric.\\
Conversely, suppose $c^k\omega$ is proportional to the metric $g^{p-k}$, then 
$\alpha_{ik}\omega_{p-i}=0$ for $k\leq i\leq p-1$. If $\alpha_{ik}=0$ then  $2p>n$ and $i<2p-n$
so that $\omega_{p-i}=0$. In  case $\alpha_{ik}\not=0$, we have $\omega_{p-i}=0$.\\
So in both cases, we have $\omega_{p-i}=0$ for $k\leq i\leq p-1$, that is 
$\omega_r=0$  for $1\leq r\leq p-k$. \ppp

\section{ Generalized Einstein Metrics}
Recall that a Riemannian manifold  is said to be Einstein if the first  contraction
 of its Riemann curvature 
tensor $R$ is proportional to the metric, that is  $cR=\lambda g$. Einstein metrics
are also known to be the critical metrics of the total scalar curvature functional once restricted
to metrics of unit volume. The gradient of the above functional is the Einstein tensor.\\
In \cite{Labbivariation}, we considered the critical metrics of the total Gauss-Bonnet
curvature functionals once restricted to metrics with unit volume. The resulting metrics are called $(2k)$-Einstein and are
 characterized by the condition that the  contraction of order $(2k-1)$ of Thorpe's
tensor $R^k$ is proportional to the metric, that is 
\begin{equation}\label{2ke}
c^{2k-1}R^k=\lambda g.
\end{equation}
\noindent
Special classes of the previous metrics are defined by:
\begin{definition} 
Let $0<p<2q<n$, we say that a  Riemannian $n$-manifold
 is $(p,q)$-Einstein if the $p$-th  contraction
 of  Thorpe's tensor of order $q$ 
is proportional to the metric
 $g^{2k-p}$, that is
\begin{equation}\label{hypere}
c^pR^q=\lambda g^{2q-p}.
\end{equation}
\end{definition}
We recover the usual Einstein manifolds for $p=q=1$ and the previous $(2q)$-Einstein condition
for $p=2q-1$.\\
Applying $(2q-p-1)$ contractions to the previous  equation (\ref{hypere}) we get 
$$c^{2q-1}R^q=\frac{(2q-p)!(n-1)!}{(n-2q+p)!}\lambda g,$$ 
that is the
 $(2q)$-Einstein condition above. Therefore we have for all $p$ with $0<p<2q$:\\
\begin{center}
$(p,q)$-Einstein $\Rightarrow$ $(2q)$-Einstein.
\end{center}
\par\smallskip\noindent
In particular, the  $(p,q)$-Einstein metrics are all critical metrics for the total Gauss-Bonnet
curvature functional of order $2q$ once restricted to metrics with unit volume. Also, 
 Schur's theorem in \cite{Labbivariation} 
implies that
 the function $\lambda$ in (\ref{hypere}) is then  a constant.\\
It is evident that for all $p\geq 1$, $(p,q)$-Einstein implies $(p+1,q)$-Einstein. In particular, 
the metrics with constant $q$-sectional curvature (that is the
sectional curvature of $R^q$ is constant) are $(p,q)$-Einstein for all $p$.\\
On the other hand, the $(p,q)$-Einstein condition neither implies nor is
 implied by the $(p,q+1)$-condition  as shown by the following examples:\\
  Let M be a 3-dimensional non-Einstein Riemannian manifold and $T^k$ be the $k$-dimensional flat torus, $k\geq 1$,
 then the Riemann curvature tensor $R$ of the Riemannian product $N=M\times T^k$ satisfies $R^q=0$ for $q\geq 2$.
In particular, $N$ is $(p,q)$-Einstein for all $p\geq 0$ and $q\geq 2$ but it
 is not $(1,1)$-Einstein.\\
On the other hand, let  $M$ be a 4-dimensional Ricci-flat but not flat manifold
(for example a K3 surface endowed with the Calabi-Yau metric), then the Riemannian product $N=M\times T^k$ is 
$(1,1)$-Einstein
but not $(q,2)$-Einstein for any $q$ with $0\leq q\leq 3$.\\
\par\medskip\noindent
It results directly from theorem 5.6 of \cite{Labbidoubleforms} and the previous proposition
 the following:
\begin{theorem}\label{thirdproposition}
 Let $0<p<2q$. The following statements are equivalent for a Riemannian
  manifold
 $(M,g)$ of dimension
$n\geq p+2q$:
\begin{enumerate}
\item The manifold $(M,g)$ is  $(p,q)$-Einstein.
\item The divergence free tensor $R_{(p,q)}=*\frac{g^{n-2q-p}}{(n-2q-p)!}R^q$ is
 proportional
 to the metric $g^p$.
\item The components $\omega_r$ of $R^q$ (with respect to the orthogonal splitting
(\ref{ort:decom}) for Thorpe's tensor $R^q$) vanish for $1\leq r\leq 2q-p$.
\end{enumerate}
\end{theorem}
The divergence free tensor $R_{(p,q)}$ used in the second part of the previous proposition
(which is nothing but the $(p,q)$-curvature tensor \cite{Labbivariation})
 generalizes the Einstein tensor obtained for $p=q=1$. \\
It results from the second part of the above proposition that
\begin{corollary}
For a 
$(p+2q)$-dimensional manifold, being  $(p,q)$-Einstein is equivalent to constant
$q$-sectional curvature.
\end{corollary}
Recall that constant $q$-sectional curvature means that the sectional curvature of Thorpe's tensor $R_{2q}$ is constant.\\
  The previous result is the analogous of the fact that in dimension $3$, Einstein manifolds
are those with constant sectional curvature.\\
The last part of the previous theorem generalizes a well known characterization of
 Einstein manifolds ($p=q=1$) by the vanishing of the component $\omega_1$ in the orthogonal decomposition
for the Riemann curvature tensor $R$.
\par\medskip\noindent
Using generalized Avez-type formulas we proved in \cite{LabbiYamabe} the following
\begin{theorem}[\cite{LabbiYamabe}] Let $(M,g)$ be a Riemannian  manifold of even dimension $n=2k+2$. Suppose the metric $g$ is $(2k-2,k)$-Einstein then the Gauss-Bonnet integrand $h_n$ of $(M,g)$ equals
\[h_n=\frac{h_{n-2}}{n(n-1)}Scal.\]
In particular,
\[ \chi(M)=ch_{n-2}\int_M Scal\,{\rm dvol}.\]
Where $c$ is a positive constant, $h_{n-2}$ is a constant : the $(2k-2,k)$-Einstein constant. Precisely it is determined by the condition $c^{2k-2}R^k=\frac{(n-2)!h_{n-2}}{2n(n-1)}g^2$, and $Scal$ denotes as usual the scalar curvature.
\end{theorem}
The previous theorem shows that for a $(2k-2,k)$-Einstein manifold of dimension $2k+2$ with $h_{n-2}\neq 0$, the integral of the scalar curvature is a topological invariant like in the two dimensional case.
\par\medskip\noindent
Using the same formulas one can prove  the following:
\begin{theorem} Let $(M,g)$ be a conformally flat $2k$-Einstein manifold of even dimension dimension $n=2k+2\geq 4$. Then the Gauss-Bonnet integrand $h_n$ of $(M,g)$ equals
\[h_n=\frac{h_{n-2}}{n(n-1)}Scal.\]
In particular,
\[ \chi(M)=ch_{n-2}\int_M Scal\,{\rm dvol}.\]
Where $c$ is a positive constant, $h_{n-2}$ is a constant : the $2k$-th Einstein constant. Precisely it is determined by the condition $c^{2k-1}R^k=\frac{(n-2)(n-3)!h_{n-2}}{n}g$, and $Scal$ is the scalar curvature.
\end{theorem}
Proof. Recall that for a conformally flat metric we have $R=Ag$ where $A$ is the Schouten tensor. On the other hand, for a $2k$ Einstein metric we have
$c^{2k-1}R^k=\lambda g$ with $\lambda=\frac{(n-2)h_{n-2}}{n}$. These two facts together with the generalized Avez type formula of \cite{LabbiYamabe} show that
\begin{equation*}
\begin{split}
h_{2k+2}=&(2k-1)\langle \lambda g,A\rangle-\langle \lambda g,cR\rangle +h_{2k}h_2\\
=& (2k-1)\lambda cA- \lambda c^2R +h_{2k}h_2\\
=& \frac{(n-2k)(n-2k-1)}{n(n-1)}h_{2k}h_2.
\end{split}
\end{equation*}
This completes the proof of the theorem. \ppp
\par\medskip\noindent
We bring the attention of the reader to a recent work of Lima and Santos \cite{LimaSantos} where they prove interesting results about the moduli space of $2k$ Einstein structures.

\subsection{Hyper $2k$-Einstein metrics}
For $q$ fixed, the $(1,q)$-Einstein  condition is the strongest among all the 
other $(p,q)$ conditions. In the rest of this section we shall emphasize  some
properties of this particular condition.\\
\begin{definition} 
Let $0<2q<n$, we shall say that a  Riemannian $n$-manifold
 is hyper $(2q)$-Einstein if the first contraction
 of Thorpe's tensor  $R^q$
is proportional to the metric
 $g^{2q-1}$, that is a $(1,q)$-Einstein manifold.
\end{definition}
The next proposition provides  topological and geometrical obstructions
 to the existence of  hyper
 $(2q)$-Einstein metrics: 
\begin{theorem}\label{thm} Let  $k\geq 1$ and  $(M,g)$ be a hyper  $(2k)$-Einstein  manifold 
of dimension
   $n\geq 4k$. 
Then  the Gauss-Bonnet curvature $h_{4k}$ of  $(M,g)$ is nonnegative. Furthermore, 
$h_{4k}\equiv 0$ if and only if  $(M,g)$ is $k$-flat.\\
In particular, a compact hyper $(2k)$-Einstein  manifold of dimension   $n= 4k$  has its
 Euler-Poincar\'e characteristic nonnegative. Furthermore,
it is zero if and only if the metric is  $k$-flat.
\end{theorem}
Recall that $k$-flat means that the sectional curvature of  $R^k$ is identically zero.
The previous proposition generalizes a well known obstruction to the existence of Einstein metrics
 in dimension four due to  Berger. \\
\par\medskip\noindent
{\sl Proof.} For a hyper $(2k)$-Einstein metric, the components $R^q_r$ vanish
for $1\leq r\leq 2k-1$ by the previous proposition. Next, for $n\geq 4k$,
 the generalized Avez formula, 
see corollary 6.5 in
\cite{Labbidoubleforms}, shows that
\begin{equation}
h_{4k}=\frac{1}{(n-4k)!}\left\{ n! ||R^k_0 ||^2 +(n-2k)! ||R^k_{2k} ||^2
\right\}.
\end{equation}
Where $R^k_i$ denotes the component of $R^k$ with respect to the orthogonal decomposition
(\ref{ort:decom}). In particular, $h_{4k}\geq 0$ and it vanishes if and only if
 $R^k_i=0$ for all $0\leq i\leq 2k$, that is $R^k\equiv 0$.\\
If $n=4k$, then $h_n$ is up to a positive constant the Gauss-Bonnet integrand. 
This completes the proof. \ppp
\par\medskip\noindent
A double form $\omega$ is said to be harmonic if it is closed and co-closed that is
 $D\omega=\delta \omega=0.$ Where $D$ is the second Bianchi map (or equivalently, the vector
valued exterior derivative)
and $\delta=c\tilde D+\tilde D c$ its "formal adjoint", $\tilde D$ being the 
adjoint second Bianchi map, see \cite{Kulk,Labbivariation}. 
\begin{proposition} Let $k\geq 1$, then for a hyper  $(2k)$-Einstein manifold, Thorpe's curvature tensor $R^k$ is a harmonic double form.
\end{proposition}
{\sl Proof.} Since $R$ satisfies the second Bianchi identity then $DR^k=0$. On the other hand,
$$\delta R^k=(c\tilde D+\tilde D c)R^k=\tilde D(cR^k)=\tilde D(\lambda g^{2k-1})=0.$$
\ppp
\\
The converse in the previous proposition is not generally true. One needs to impose 
extra conditions
on the manifold, similar to the ones in \cite{Bourguignon},
 in order to ensure the generalized Einstein condition.
\par\medskip\noindent
For $n\geq 4k$, following \cite{Kulk,Nasu}, an $n$-dimensional Riemannian manifold is said
 to be $k$-conformally flat if Thorpe's tensor $R^k$ is divisible by the metric $g$. Equivalently, the
component $(R^k)_{2k}$ of $R^k$ with respect to the orthogonal splitting
(\ref{ort:decom}) vanishes. The proof of the following proposition is straightforward 
if one uses theorem \ref{thirdproposition}
\begin{proposition} A hyper $(2k)$-Einstein manifold which is $k$-conformally flat has
constant $k$-sectional curvature.
\end{proposition}
The previous proposition generalizes a similar result about usual Einstein conformally flat
manifolds obtained for $k=1$.\\
\section{Thorpe's Condition}
In this section we consider different generalizations of the Einstein condition suggested by the 
work of Thorpe \cite{Thorpe}.\\
Recall that in dimension four, Einstein metrics are characterized by the invariance of their
curvature tensor under
the Hodge star operator in the sens that $*R=R$, where $R$ is seen as a double form and 
$*$ is
the generalized Hodge star operator \cite{Labbidoubleforms}. In higher even dimensions
 we proved a similar characterization:
 \begin{proposition}[\cite{Labbidoubleforms}]
A Riemannian manifold $(M,g)$ of dimension $n=2p\geq 4$ 
is Einstein
if and only if its Riemann curvature tensor $R$ satisfies
\begin{equation} *g^{p-2}R=g^{p-2}R.
\end{equation}
\end{proposition}
For $(4k)$-dimensions, Thorpe considered the following similar condition but on the generalized
 Thorpe's  tensor instead, namely
$*R^k=R^k.$\\
 For a $(4k)$-dimensional
compact and oriented Riemannian manifold, he proved that the previous condition implies 
\begin{equation}
\chi\geq \frac{(k!)^2}{(2k)!}|p_k|.
\end{equation}
Furthermore, $\chi=0$ if and only if the manifold is $k$-flat.
Where $p_k$ denotes the $k^{th}$ Pontrjagin number of the manifold.\\
If $R^k=\sum_{i=0}^{2k}g^i\omega_{2k-i}$ is the decomposition of $R^k$ following the orthogonal
splitting (\ref{ort:decom}), then Thorpe's condition $*R^k=R^k$ is equivalent  
 \cite{Labbidoubleforms} to
the vanishing
of $\omega_r$ for $r$ odd and $1\leq r\leq 2k-1$. This is in turn can be seen to be
 equivalent to the condition that
$c^rR^k$ is divisible per the metric (that is $c^rR^k=g {\bar{\omega}}$) for $r$ 
odd and $1\leq r\leq 2k-1$. These restrictions on $R^q$ are guaranteed  
in arbitrary even  dimension $n=2p\geq 4q$ if one assumes  the following generalized
 Thorpe type condition,
see theorem 5.8 in \cite{Labbidoubleforms} 
\begin{equation} *g^{p-2q}R^q=g^{p-2q}R^q.
\end{equation}
The previous generalization of Thorpe's condition forces the Gauss-Bonnet curvature
$h_{4k}$ to be nonnegative in dimensions $n\geq 4k$ (theorem 6.7 of \cite{Labbidoubleforms}).\\
\par\medskip\noindent
It results from the previous discussion that for a $(4k)$-dimensional manifold,
the hyper $(2k)$-Einstein condition is stronger than
Thorpe's condition $*R^k=R^k$. In higher arbitrary even dimensions, the $(p,q)$-Einstein condition
implies the generalized Thorpe's condition $*g^{p-2q}R^q=g^{p-2q}R^q$. In particular, theorem
\ref{thm} is then a consequence of the above result of Thorpe in the special case where  $n=4k$. Also, 
 theorem 6.7 of \cite{Labbidoubleforms} implies the conclusions of theorem \ref{thm}
for  even dimensions  higher than $4k$.\\
 On the other hand
Thorpe's condition and its generalization are actually 
complicated expressions once written in terms of
the contractions of $R^q$ as it is shown below.\\
The proof of the following proposition is a direct application of formula (15) 
of \cite{Labbidoubleforms} 
\begin{proposition}
 A Riemannian manifold of dimension $n=2p\geq 4q$ satisfies the generalized
Thorpe's condition (that is $*g^{p-2q}R^q=g^{p-2q}R^q$) if and only if
\begin{equation}
\sum_{r=1}^{2q}\frac{(-1)^r(p-2q+1)!}{r!(p-2q+r)!}g^{r-1}c^rR^q=0.
\end{equation}
In particular, if $n=4q$, Thorpe's condition is equivalent to
\begin{equation}
\sum_{r=1}^{2q}\frac{(-1)^r}{(r!)^2}g^{r-1}c^rR^q=0.
\end{equation}
\end{proposition}

\subsection*{Final remarks and questions}
Our study of $(p,q)$ Einstein manifolds is somehow premature, in the sens that still we do not fully understand the model spaces where the $q$-th Thorpe's sectional curvature is constant (that are the $(0,q)$ Einstein metrics).\\
It is an interesting question to classify the Riemannian manifolds of constant $k$-th Thorpe's sectional curvature.\\
It would also be interesting to generalize Hamilton's Ricci flow techniques to the various generalized Ricci curvatures introduced here in this paper.

\subsection*{Acknowledgment}
I would like to thank the referee for useful comments about the first version of this paper.

\vspace{2cm}
\noindent
Labbi M.-L.\\
  Department of Mathematics,\\
 College of Science, University of Bahrain,\\
  P. O. Box 32038 Bahrain.\\
  E-mail: labbi@sci.uob.bh


\begin{thebibliography}{9}

\bibitem{Bourguignon} Bourguignon, J. P. \emph{Les vari\'et\'es de dimension $4$ \`a 
signature non nulle dont la courbure est harmonique sont d'Einstein}, Invent. Math. 63, 263-286
(1981).
\bibitem{Kulk} Kulkarni, R. S., \emph{On Bianchi Identities},
 Math. Ann. 199, 175-204(1972).
\bibitem{Labbidoubleforms} Labbi, M.-L., \emph{ Double forms, curvature structures and the
$(p,q)$-curvatures}, Transactions of the American Mathematical Society, 357, 
n10, 3971-3992 (2005).
\bibitem{Labbivariation} Labbi, M. L. \emph{ Variational properties of the 
Gauss-Bonnet curvatures}, Calc. var. (2008) 32: 175-189.
\bibitem{LabbiYamabe} Labbi, M. L. \emph{About the $h_{2k}$ Yamabe problem}, preprint arXiv:0807.2058v1 [math.DG].
\bibitem{LimaSantos} Lima, L. L. and Santos, N. L., \emph{Infinitesimal deformations of $2k$-Einstein structures}, preprint.
\bibitem{Nasu} Nasu, T., \emph{On conformal invariants of higher order}, 
Hiroshima Math. J. 5 (1975), 43-60.
\bibitem{Thorpe} Thorpe, J. A., \emph{Some remarks on the Gauss-Bonnet
integral}, Journal of Mathematics and Mechanics, Vol. 18, No. 8 (1969).

\end{thebibliography}
\end{document}